\documentclass{amsart}

\usepackage{latexsym}
\usepackage{euscript}

\newenvironment{pf}{\begin{trivlist} \item[] {\it Proof.} \ }{\qed \end{trivlist} } 
\newenvironment{pf*}[1]{\begin{trivlist} \item[] {\it #1.} \ }% {\qed \end{trivlist} }

\newtheorem{theorem}{Theorem}%[section]

\newtheorem{proposition}[theorem]{Proposition}

\newtheorem{assertion}[theorem]{Assertion}

\newtheorem{definition}[theorem]{Definition}

\theoremstyle{remark}

\newtheorem{remark}[theorem]{Remark}

\newcommand{\w}{\wedge}
\newcommand{\om}{\omega}

\newcommand{\lam}{\lambda}

\newcommand{\calI}{{\mathcal I}}
\newcommand{\calJ}{{\mathcal J}}
\newcommand{\calK}{{\mathcal K}}

\begin{document}

\title[On the intermediate integral for Monge-Amp\`ere equations]{On the 
intermediate integral for Monge-Amp\`ere equations} 

\author {Jeanne Nielsen Clelland}
\address{School of Mathematics, Institute for Advanced Study\\ Princeton, NJ 
08540}
\email{jnc@math.ias.edu}

\thanks{This research was supported in part by NSF grant DMS-9627403.}

\begin{abstract}
Goursat showed that in the presence of an intermediate integral, the 
problem of solving a second-order Monge-Amp\`ere equation can be 
reduced to solving a first-order equation, in the sense that the 
generic solution of the first-order equation will also be a solution 
of the original equation.  An attempt by Hermann to 
give a rigorous proof of this fact contains an error; we show that 
there exists an essentially unique counterexample to Hermann's 
assertion and state and prove a correct theorem.

\end{abstract}

\maketitle

The purpose of this paper is to correct an error in the literature 
regarding Goursat's method of the intermediate integral for 
Monge-Amp\`ere equations.  Recall that a {\em Monge-Amp\`ere 
equation} is a second-order PDE of the form
\begin{equation}
A(z_{xx} z_{yy} - z_{xy}^2) + B z_{xx} + 2C z_{xy} + D
z_{yy} + E = 0  \label{MAeq}
\end{equation}
where the coefficients $A,B,C,D,E$ are functions of the five variables 
$x,y,z,z_x,z_y$.  Such an equation can be described geometrically as 
follows: let $M = \mathbb{R}^5$ with coordinates $x,y,z,p,q$, and 
consider the differential ideal $\calI$ on $M$ generated by the forms
\begin{gather*}
\theta = dz - p\,dx - q\,dy \\
d\theta = -dp \w dx - dq \w dy \\
\Omega = A\, dp \w dq + B\, dp \w dy + C\, (dx \w dp + dq \w dy) + D\, 
dx \w dq + E\, dx \w dy .
\end{gather*}
The coefficients $A,B,C,D,E$ are now thought of as functions of 
the variables $x,y,z,p,q$.  2-dimensional submanifolds $i:N^2 
\hookrightarrow M$ satisfying $i^*(\theta) = i^*(d\theta) = 
i^*(\Omega) = 0$,\ $i^*(dx \w dy) \neq 0$ are in one-to-one correspondence 
with solutions of \eqref{MAeq}.

This motivates the following geometric definition.

\begin{definition}
A {\em Monge-Amp\`ere system} on a 5-manifold $M^5$is a differential 
ideal $\calI \subset \Omega^*(M)$  which is generated by a contact 1-form 
$\theta$ (i.e., a 1-form $\theta$ with the property that $\theta \w 
(d\theta)^2$ never vanishes on $M$), its exterior derivative $d\theta$, 
and another 2-form $\Omega$ with the property that $\Omega$ and 
$d\theta$ are linearly independent modulo $\theta$ at each point of 
$M$.  
\end{definition}

This definition includes systems that arise from Monge-Amp\`ere 
equations as described above, as well as more general systems.  An {\em 
integral manifold} of a Monge-Amp\`ere system is a 2-dimensional 
submanifold $i:N^2 \hookrightarrow M$ such that $i^*(\theta) = 
i^*(d\theta) = i^*(\Omega) = 0$.  If the system arises from a PDE as 
above, then solutions of the PDE are typically in one-to-one 
correspondence with those 
integral manifolds of $\calI$ that satisfy a certain independence condition.  
However, in some contexts more general integral manifolds may be of interest 
as well.
For any Monge-Amp\`ere system, there always exist local coordinates 
$x,y,z,p,q$ on $M$ in which $\calI$ takes the form described 
above, although there may not be a global coordinate expression of 
this form.  In terms of the PDE \eqref{MAeq}, linear independence 
of $\Omega$ and $d\theta$ modulo $\theta$ means that the coefficients 
$A,B,C,D,E$ never vanish simultaneously on $M$.  Further discussion of 
Monge-Amp\`ere systems may be found in \cite{BG2}.

The 2-form $\Omega$ may be replaced by any expression of the form 
$\Omega + \lam d\theta$ without changing the ideal $\calI$.  
Suppose that $\lam$ is chosen so that the resulting form $\Omega$ is 
decomposable; i.e., that there exist two 1-forms $\om^1, \om^2$ such that 
\[ \Omega = \om^1 \w \om^2. \]
There is a unique (up to nonzero multiples) such choice for $\Omega$ in 
the case that $\calI$ 
represents a parabolic equation, and there are exactly two distinct 
such choices for $\Omega$ in the case that $\calI$ represents a 
hyperbolic equation (or an elliptic equation if the forms $\om^1, 
\om^2$ are allowed to assume complex values).

Having chosen $\Omega$ as above, let $\calJ$ be the differential system 
defined by 
\[ \calJ = \{\theta, \om^1, \om^2\} . \]
$\calJ$ is called a {\em characteristic system} of $\calI$, and
the integral curves of $\calJ$ are {\em characteristics} of the 
system $\calI$.  This definition agrees with the usual notion
of characteristics.

Now suppose that there exists a function $F$ on $M$ such that $dF 
\in \calJ$.  Goursat \cite{G} calls such a function a {\em first-order 
intermediate integral} for $\calI$.  The reason for this terminology 
is as follows: let $F$ be given in local coordinates as 
$F(x,y,z,p,q)$, and suppose that $N^2 \subset M$ is an integral manifold 
of the system $\calK$ defined by $\{\theta, d\theta, F \}$.  If 
$N$ satisfies the independence condition $dx \w dy 
\neq 0$, then it represents a solution of the first-order PDE
\begin{equation}
F(x, y, z, z_x, z_y) = 0. \label{PDE1}
\end{equation}  
Goursat shows that generically, any 
such $N$ is in fact an integral manifold of $\calI$.  Goursat does 
allow for the possibility of exceptional $N$ that do not have
this property (``Le raisonnement ne pourrait \^etre en d\'efaut que 
si, pour l'int\'egrale consid\'er\'ee, un des facteurs $\lam, \nu, 
\mu$ devenait ind\'etermin\'e,'' \cite{G} p. 60).  

From the point of view of PDEs, this result gives a method 
for finding solutions to equation \eqref{MAeq}: first find an intermediate 
integral $F$, and then solve the 
first-order PDE \eqref{PDE1}.  Any solution of \eqref{PDE1}
should be a solution of the original equation \eqref{MAeq}.  Goursat 
uses this method to show, for example, that the solutions of the equation
\[ z_{xx}z_{yy} - z_{xy}^2 = 0 \]
are exactly the developable surfaces.

In \cite{H}, Hermann attempts to give a rigorous proof of Goursat's 
result for an arbitrary integral manifold of $\calK$, but the proof 
contains an error.  In the remainder of the paper we will discuss 
Hermann's argument and the error it contains, describe
counterexamples, and state and prove a correct theorem.

Theorem 8.1 of \cite{H} may be stated in local coordinates as follows:

\begin{assertion}\label{Hermann}
Suppose that $F$ is an intermediate integral for $\calI$ such
that $dF$ never vanishes on $M$.  Let $i:N^2 \hookrightarrow M$ be an integral 
manifold of $\calK$ that satisfies the independence condition 
$i^*(dx \w dy) \neq 0$.  Then $N$ is an integral manifold of $\calI$.
\end{assertion} 

The independence condition actually does not enter into the proof, and 
it is convenient to omit this hypothesis in what follows so that we may 
consider equivalence of systems under contact transformations.  
Two Monge-Amp\`ere systems $(M, \calI),\ (\bar{M}, \bar{\calI})$ are said 
to be {\em contact equivalent} if there is a diffeomorphism
\linebreak $\psi:M 
\rightarrow \bar{M}$ such that $\psi^*\bar{\calI} = \calI.$  This is a 
natural notion of equivalence for differential systems, but 
such transformations do not necessarily preserve independence conditions.  

The argument given in \cite{H} is as follows: since $F$ is an intermediate 
integral of $\calI$, we have
\begin{equation}
dF = a\theta + b \om^1 + c \om^2, \label{dFeq}
\end{equation}
where by hypothesis $a,b,c$ never vanish simultaneously.  The hypotheses 
also imply that $i^*(\theta) = i^*(dF) = 0$; therefore, 
\begin{equation}
i^*(b \om^1 + c \om^2) = 0. \label{dFpullback}
\end{equation}
The only way that $N$ can fail to be an integral manifold of $\calI$ 
is if $i^*(\om^1 \w \om^2) \neq 0$, and by \eqref{dFpullback} this can 
only happen if $b = c = 0$.  But taking the exterior derivative of 
both sides of \eqref{dFeq} shows that this impossible unless we also 
have $a=0$, which is a contradiction.

The error is in this last step.  The conclusion would be true if we knew 
that $b=c=0$ held on all of $M^5$, but in fact we only know that 
$i^*(b) = i^*(c) = 0$, i.e., that $b$ and $c$ vanish on a 
2-dimensional submanifold of $M$.  This is perfectly possible; a 
counterexample (admittedly a very degenerate one) to the assertion is as 
follows: take $M = \mathbb{R}^5,\ \om^1 = dx, \ \om^2 = dy, \ F = z$, and 
let $i:N^2 \hookrightarrow M$ be the surface defined by $z = p = q = 0$.
Then 
\[ dF = \theta + p\om^1 + q\om^2 \] 
and $N$ is an integral manifold 
of  $\calK$, but $i^*(\om^1 \wedge \om^2) = dx 
\wedge dy \neq 0$, so $N$ is not an integral manifold of $\calI$. 

\begin{remark}
As described above, $\calI$ represents the equation $1=0$, 
which is not exactly an interesting partial differential equation.  But 
in the coordinate system $(X,Y,Z,P,Q)$ defined by
\[ X = p,\ \ \ Y = q, \ \ \ Z = z - p\,dx - q\,dy, \ \ \ P = 
-x, \ \ \ Q = -y, \]
$\calI$ represents the equation
\[ Z_{XX}Z_{YY} - Z_{XY}^2 = 0. \]
The exceptional integral manifold $N$ in these coordinates is not a 
solution of the equation in the usual sense; it consists of 
all planes passing through the point $(0,0,0) \in \mathbb{R}^3$.
\end{remark}

The good news, however, is that up to local contact equivalence this is 
the only counterexample.  We will demonstrate this in two steps.  
Note that this example has the property that 
$\calJ$ defines an integrable distribution.  First we will show that any 
counterexample for which $\calJ$ is integrable is locally contact 
equivalent to the example above.  Next we will show that 
Assertion \ref{Hermann} holds if $\calJ$ is not integrable.

\begin{proposition}
Suppose that $\calJ$ is integrable, and let 
$F$ be an intermediate integral for $\calI$ such that $dF$ never 
vanishes on $M$.  If $N^2$ is an integral manifold of $\calK$ 
but not of $\calI$, then there exist local coordinates $x,y,z,p,q$ 
such that $\calI$ is generated by the forms 
\begin{gather*}
\theta = dz - p\, dx - q\, dy \\
d\theta = -dp \w dx - dq \w dy \\
\Omega = dx \w dy,
\end{gather*}
$F = z$, and $N$ is defined by the equations $z = p = q = 0$.
\end{proposition}

\begin{pf}
Since $\calJ$ is integrable, the Frobenius theorem implies that 
locally there exist independent functions $X,Y,Z$ on $M$ such that
\[ \calJ = \{dX, \, dY, \, dZ \}. \]
Moreover, the Pfaff theorem implies that these functions can be chosen 
so that up to a nonzero multiple,
\[ \theta = dZ - P\, dX - Q\, dY \]
for some functions $P,Q$ on $M$.  Since $\theta \w (d\theta)^2 \neq 
0$, the functions $X,Y,Z,P,Q$ are independent and so form a local 
coordinate system on $M$.  (See \cite{BCG3} for details.)  Since 
$\Omega = \om^1 \w \om^2$ with $\om^1,\om^2 \in \calJ$, we must have
$\Omega \equiv \lam\,dX \w dY \mod{\theta}$ for some nonvanishing $\lam$.

Now any function $F(X,Y,Z)$ is an intermediate integral for $\calI$; 
conversely, every intermediate integral must have this form.  For any 
such function, we can write
\[ dF = F_Z\, \theta + (F_X + P F_Z)\, dX + (F_Y + Q F_Z)\, dY. \]
If $N^2$ is an integral manifold of $\calK$ but not of $\calI$, then 
$N$ must be defined by the equations
\[ F = F_X + P F_Z = F_Y + Q F_Z = 0. \]
Since $dF$ never vanishes, it must also be true that $F_Z \neq 0$ on 
$N$.  In a neighborhood of $N$, define new coordinates $(x,y,z,p,q)$ by
\begin{gather*}
x = X \\
y = Y \\
z = F(X,Y,Z) \\
p = F_X + P F_Z \\
q = F_Y + Q F_Z.
\end{gather*}
We have
\[ dz - p\, dx - q\, dy = F_Z\, \theta, \]
so this change of variables is a contact transformation.  
Furthermore, $F = z$ and
$N$ is defined by the equations $z = p = q = 0.$
\end{pf}

\begin{theorem}  
Suppose that the system $\calJ$ is not integrable and that either
\begin{enumerate}
\item{$d\theta \equiv 0 \mod{\calJ}$, or}
\item{$\theta \w d\theta \w \Omega$ never vanishes.}
\end{enumerate}
Let $F$ be an intermediate integral for $\calI$ such
that $dF$ never vanishes on $M$.  If $i:N^2 \hookrightarrow M$ is any integral 
manifold of $\calK$, then $N$ is an integral manifold of $\calI$.  
\end{theorem}

The assumptions on $d\theta$ are essentially constant type 
assumptions; the first assumption holds for parabolic equations, and 
the second holds for hyperbolic or elliptic equations.

\begin{pf}
Choose 1-forms $\om^3, \om^4$ on $M$ so that $\{\theta, \om^1, \om^2, 
\om^3, \om^4 \}$ is a basis for the cotangent space at each point of 
$M$.  Since $F$ is an intermediate integral for $\calI$, $dF$ must 
lie in the last derived system of $\calJ$.  (See \cite{BCG3} for a 
discussion of derived systems.)

First suppose that assumption 1 holds.  Then
\[ \left. \begin{array}{l} 
d\theta \equiv 0 \\
d\om^1 \equiv r_1\, \om^3 \w \om^4 \\
d\om^2 \equiv r_2\, \om^3 \w \om^4 
\end{array} \right\}  \mod{\calJ}  \]
where by the nonintegrability hypothesis $r_1, r_2$ never vanish 
simultaneously.  By choosing a different 
basis for $\{\om^1, \om^2\}$ if necessary, we can assume that $r_1 
\equiv 0, r_2 \neq 0$, so that
\[ \calJ^{(1)} = \{\theta, \om^1\}. \]
Now since $\theta$ is a contact form, $d\theta$ has no linear 
divisors modulo $\theta$.  Therefore the expression $d\theta 
\mod{\calJ^{(1)}}$ cannot vanish at any point of $M$.  So if $\calJ^{(2)} = 
(\calJ^{(1)})^{(1)}$ is nonzero, then it has dimension 1 and can be 
expressed in the form
\[ \calJ^{(2)} =  \{ \om^1 - s\theta \} \]
for some function $s$.  Therefore, any intermediate integral $F$ 
satisfies
\[ dF = a(\om^1 - s\theta) \]
with $a \neq 0$, and the argument of \cite{H} is valid.

Next suppose that assumption 2 holds.  Then
\[ \left. \begin{array}{l} 
d\theta \equiv r_0\, \om^3 \w \om^4 \\
d\om^1 \equiv r_1\, \om^3 \w \om^4 \\
d\om^2 \equiv r_2\, \om^3 \w \om^4 
\end{array} \right\}  \mod{\calJ}  \]
where $r_0$ never vanishes.  So $\calJ^{(1)}$ is 2-dimensional 
and can be expressed in the form
\[ \calJ^{(1)} =  \{ \om^1 - s_1\theta,\ \om^2 - s_2 \theta \} \]
for some functions $s_1, s_2$.  Therefore, any intermediate integral $F$ 
satisfies
\[ dF = a(\om^1 - s_1\theta) + b(\om^2 - s_2\theta) \]
where $a,b$ never vanish simultaneously, and the argument of \cite{H} is 
valid.

\end{pf}

\end{document}